\documentclass[letterpaper,10pt]{article}
\usepackage{enumitem}
\usepackage{amsmath,amssymb} 
\usepackage{cases}
\usepackage{comment}
\usepackage{hyperref}
\usepackage{longtable}
\usepackage{amsfonts}
\usepackage{graphicx}
\graphicspath{{/}}
\usepackage[caption=false]{subfig}
\usepackage{listings}
\usepackage{anysize}
\marginsize{0.85in}{0.85in}{0.6in}{0.6in} 
\usepackage{booktabs}
\usepackage{url}
\usepackage{caption}
\newcommand{\besselroot}{j_\text{\tiny 0\hspace{-0.2ex}1}}
\newcommand{\lambdaodot}{\rho}
\usepackage{pxfonts}

\newcommand\eccentricity{\mathrm{e}} 
\newcommand\e{\eccentricity}

\usepackage{wrapfig}
\begin{document}
\title{The fundamental Laplacian eigenvalue of the\\ ellipse with
  Dirichlet boundary conditions} \author{Robert Stephen Jones
  \thanks{\href{mailto:rsjones7@yahoo.com}{rsjones7@yahoo.com} ,
    {\small Independent Researcher, Sunbury, Ohio}}} 

\maketitle {\abstract{\raggedright

    In this project, I examine the lowest Dirichlet eigenvalue of the
    Laplacian within the ellipse as a function of eccentricity.  Two
    existing analytic expansions of the eigenvalue are extended: Close
    to the circle (eccentricity near zero) nine terms are added to the
    Maclaurin series; and near the infinite strip (eccentricity near
    unity) four terms are added to the asymptotic expansion. In the
    past, other methods, such as boundary variation techniques, have
    been used to work on this problem, but I use a different
    approach~--~which not only offers independent confirmation of
    existing results, but may be used to extend them. My starting
    point is a high precision computation of the eigenvalue for
    selected values of eccentricity. These data are then fit to
    polynomials in appropriate parameters yielding high-precision
    coefficients that are fed into an LLL integer-relation algorithm
    with forms guided by prior results.  }} \raggedright
\section*{Introduction}
Let $\Omega$ be an elliptical region with boundary
$\partial\Omega$, as shown in Fig.~\ref{fig:ellipse}, and where
some elementary but relevant properties are given in the Appendix.

The Dirichlet Laplacian eigenvalue problem within this region is
defined by
  \begin{equation}
  \Delta \psi + \lambda \psi=0 \quad \text{in~$\Omega$\qquad with}\quad \psi =0 \text{\quad on~$\partial\Omega$}
  \end{equation}
where,~in~general,~there~exists a non-accumulating, infinite tower of
real eigenvalues
\begin{equation}
  0 < \lambda_0 < \lambda_1 \le \lambda_2 \le \lambda_3 \le \cdots
\end{equation}
and corresponding eigenfunctions, $\psi_n\in L^2$, orthonormalizable
as $\int_\Omega\psi_m\psi_n=\delta_{mn}$.  This classic problem has a
long history with many results, but only relevant techniques and
results are reported here.

This project is limited to examining how the lowest (fundamental)
eigenvalue, $\lambda_0$, behaves as a function of ellipse
parameters. Specifically, it confirms and extends two series: A
Maclaurin series for small values of eccentricity ($\e\approx 0$) near the
circle, and an asymptotic series for large values of eccentricity ($\e\lesssim 1$) as
the highly-elongated ellipse degenerates into the \textit{infinite strip}.

\begin{figure}
\centering
  \includegraphics[width=.5\linewidth]{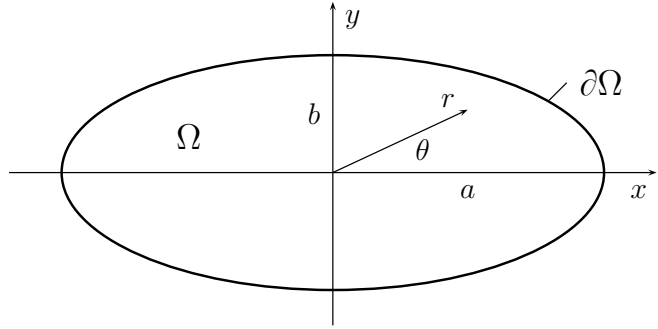}
  \caption{The ellipse.}\label{fig:ellipse}
\end{figure}

The (un-normalized) eigen-solutions for these two extremes are very
well known and elementary,
\begin{subnumcases}{\left\{ \lambda_0, \psi_0\right\}=}
  \left\{ \left(\frac{\besselroot}{R}\right)^2,\ J_0\left(\frac{\besselroot r}{R}\right)\right\} &
\text{if\quad$\eccentricity=0$\qquad $R$-radius circle,\quad $r=\sqrt{x^2+y^2}\le R$} 
 \label{eq:limitcircle}
  \\ \left\{\displaystyle \left(\frac{\pi}{2\varepsilon}\right)^2,\,
  \cos\left(\frac{\pi y}{2\varepsilon} \right)\right\} & \text{if\quad
    $\eccentricity= 1$\qquad infinite strip,\quad
    $|y|\le\varepsilon$\qquad (this $\psi_0\notin L^2$)\label{eq:limitstrip}
 } 
    \label{eq:limits}
\end{subnumcases}
where $J_0$ is the Bessel function of the first kind of order zero, and
$\besselroot\approx 2.4048$ its first root.

As noted, in the limit $\e\to 1^{-}$, the eigenfunction loses its $L^2$
property. For us, that is not a problem: We shall be interested in the
eigenvalue, $[\pi/(2\varepsilon)]^2$, since this becomes the lowest order
term in the asymptotic expansion, Eq.~(\ref{eq:resultB}), where
$\varepsilon$ of~(\ref{eq:limitstrip}) shall correspond to the stretch 
factor of the ellipse per~Eq.~(\ref{eq:ellipse}).


\section*{Results}

First, consider the Maclaurin series of $\lambda_0$ in powers of eccentricity $\e$, 
\begin{multline}
A=\pi:\qquad \frac{\lambda_0}{\lambdaodot} = 
\sum_{\nu=0}^\infty C_\nu(\lambdaodot)\, \e^{2\nu} = 1 +
\left[\frac{\lambdaodot-2}{32}\right]\left(\e^4 + \e^6\right) + \left\{\frac{- 7\rho^3 + 58\rho^2 + 832\rho - 1792}{32768}\right\}\,\e^8 \\
\mbox{} + \left\{\frac{- 7\rho^3 + 58\rho^2 + 320\rho - 768 }{16384}\right\}\,\e^{10}
 + \text{[See Table~\ref{tab:resultA}]}+\cdots\qquad
  \label{eq:resultA}
\end{multline}
where $\lambdaodot=\besselroot^2\approx 5.7831$ is the fundamental
eigenvalue of the unit-radius circle per Eq.~(\ref{eq:limitcircle}),
and where the next eight terms, from $C_6(\lambdaodot)\,\e^{12}$ to
$C_{13}(\lambdaodot)\,\e^{26}$, are listed in Table~\ref{tab:resultA}.
The coefficients~--~expressed as rational polynomials in powers
of $\lambdaodot$~--~are exact, but, for reference, rounded numerical
values of the first thirty non-trivial coefficients appear in
Table~\ref{tab:resultAnumerical}.

About fifty years ago, in 1967, Joseph~\cite{joseph1967} was first to
publish up to $C_3(\lambdaodot)$ and noted that only even powers of
$\eccentricity$ appear. A decade or two went by (date unknown) when
Henry~\cite{henry_2005} picked up the problem again, corrected a minor
mistake in the Joseph result, and also reported that result to the
same order. At least another decade elapsed when, in 2014, $C_4(\lambdaodot)$ was
determined by Boady, Grinfeld, and Johnson~\cite{bgj2013}.%
\footnote{To line up their result with Eq.~(\ref{eq:resultA}), first
  realize they used yet a different area, $A''=\pi/\varepsilon$, so
  multiply Eq.~(\ref{eq:resultA}) by $\varepsilon=\sqrt{1-\e^2}$ and
  re-expand in powers of $\e$.  The coefficient multiplying $\e^8$
  becomes the same as theirs, $ \left[-7\lambdaodot^3+58 \lambdaodot^2
    + 192\lambdaodot-1792\right]/32768$, for example.  }  All of those
prior efforts relied on boundary variation techniques by
parametrically deforming the circle into the ellipse.

\begin{table}[!tb]
  \caption{\small Higher-order terms of the Maclaurin series for the
    ellipse per Eq.~(\ref{eq:resultA}).}
  \label{tab:resultA}
  \begin{tabular}{ccl}
    \toprule
    $\nu$ & \rule{3ex}{0ex} & \qquad\qquad$C_\nu(\lambdaodot)\,\e^{2\nu}$ \\
    \midrule
    \rule{0pt}{4.5ex}%
    6 && $\displaystyle\frac{+ 87\rho^5 - 1066\rho^4 - 12778\rho^3
       + 134676\rho^2 + 418176\rho - 1140480 }{28311552}\,\e^{12}$\\[3ex]
    7 && $\displaystyle\frac{+ 87\rho^5 - 1066\rho^4 - 2698\rho^3
      + 51156\rho^2 + 104832\rho - 329472 }{9437184}\,\e^{14}$ \\[3ex]
    8 && $\displaystyle\frac{\left\{
        \begin{array}{l} - 206061\rho^7 + 3371550\rho^6 + 44817952\rho^5 - 742073664\rho^4\\
          \quad \mbox{}-4882432\rho^3 + 21039022080\rho^2 + 30916214784\rho - 113359454208
        \end{array}\right\} }{3710851743744}\,\e^{16}$\\[3ex]
9 &&  $\displaystyle\frac{\left\{
    \begin{array}{l} -206061\rho^7 + 3371550\rho^6 + 4906528\rho^5 - 253044032\rho^4\\
      \quad \mbox{}+307986432\rho^3 + 5234098176\rho^2 + 5775556608\rho - 25027411968
         \end{array}\right\} }{927712935936}~\e^{18}$ \\[3ex]
10 && $\displaystyle\frac{\left\{
  \begin{array}{l} + 16700445\rho^9 - 342482130\rho^8 - 5150192834\rho^7 + 118301328148\rho^6 \\
    \quad\mbox{}- 200183585216\rho^5 - 4688313904000\rho^4 + 9860534272000\rho^3 + 80801390592000\rho^2\\
    \quad\quad\mbox{}+ 68584734720000\rho - 356640620544000
  \end{array}\right\} }{14843406974976000}\,\e^{20} $ \\[3ex]

11 && $\displaystyle\frac{\left\{
  \begin{array}{l} + 16700445\rho^9 - 342482130\rho^8 - 204728834\rho^7 + 37384128148\rho^6\\
    \quad\mbox{} - 126365422016\rho^5 - 962599369600\rho^4 + 2869625651200\rho^3 + 15282796953600\rho^2\\
    \quad\quad\mbox{} + 10059094425600\rho - 64012419072000
  \end{array}\right\} }{2968681394995200}\,\e^{22}$ \\[3ex]

12 &&
$\displaystyle\frac{\left\{\begin{array}{l}
  +120332513685 \rho^{11}
  -2968187062070 \rho^{10}
  -50217731403560 \rho^{9}
  +1484737085079984 \rho^{8}\\
  \quad \mbox{}
  -4817972151021312 \rho^{7}
  -77508820026886656 \rho^{6}
  +383914479592341504 \rho^{5}\\
  \quad \quad \mbox{}
  +1477542066905088000 \rho^{4}
  -6033467570651136000 \rho^{3}
  -23656611409035264000 \rho^{2}\\
  \quad\quad\quad\mbox{}
  -11993659465531392000 \rho
  +95949275724251136000
  \end{array}\right\} }{4924686192529637376000}\,\e^{24}$\\[3ex]
13 && 
$\displaystyle\frac{\left\{\begin{array}{l}
  +120332513685 \rho^{11}
  -2968187062070 \rho^{10}
  +572997966040 \rho^{9}
  +443153032753584 \rho^{8}\\
  \quad\mbox{}
  -2432501708504832 \rho^{7}
  -13031591176137216 \rho^{6}
  +94442136581505024 \rho^{5}\\
  \quad\quad\mbox{}
  +204288219832320000 \rho^{4}
  -1178391769251840000 \rho^{3}
  -3645518590771200000 \rho^{2}\\
  \quad\quad\quad\mbox{}
  -1389639688519680000 \rho
  +14537769049128960000
  \end{array}\right\} }{820781032088272896000}\,\e^{26}$\\[3ex]
\bottomrule
  \end{tabular}
 \end{table}

One of my contributions to this problem is the computation
up to $C_{13}(\rho)$, including independent verification of the coefficients up to $C_4(\rho)$.

\begin{table}[!htb]
  \centering
  \caption{\small Numerical values of the leading thirty coefficients of the Maclaurin
    series for the eigenvalue within the ellipse per
    Eq.~(\ref{eq:resultA}), and rounded to twenty decimal
    places. Also indicated is $D_\nu$, the number of digits in agreement
  between the numerical coefficient (via this linear regression of data)
  and the respective coefficient displayed in Eq.~(\ref{eq:resultA})
  with Table~\ref{tab:resultA}.  This truncated list of coefficients
  is based a fit using fifty coefficients and sixty 500-digit eigenvalues with eccentricity $\e=0.000001$
  to $0.000060$, spaced equally; and it incorporates the trivial
  $C_0=1$ and $C_1=0$.}
\label{tab:resultAnumerical}
\begin{tabular}{rlcccrlcc}
\toprule
$\nu$ &\qquad$C_\nu$ & $C_\nu/C_{\nu-1}$ & $D_\nu$ &\qquad \qquad & $\nu$ & \qquad$C_\nu$ & $C_\nu/C_{\nu-1}$\\[1ex]
\midrule
2 & 0.11822456134208701629 &  & 458  & & 17 & 0.05770190566258202267 & 0.97004 &   \\[1ex]
3 & 0.11822456134208701629 & 1.00000 & 446  & & 18 & 0.05607079818966375885 & 0.97173 &   \\[1ex]
4 & 0.11003095525016373549 & 0.93069 & 434  & & 19 & 0.05457135306802184693 & 0.97326 &   \\[1ex]
5 & 0.10183734915824045469 & 0.92553 & 423  & & 20 & 0.05318746732402259177 & 0.97464 &   \\[1ex]
6 & 0.09469809424285786691 & 0.92990 & 412  & & 21 & 0.05190553831652441797 & 0.97590 &   \\[1ex]
7 & 0.08861319050401597214 & 0.93574 & 402  & & 22 & 0.05071400061380009767 & 0.97704 &   \\[1ex]
8 & 0.08341794996585013471 & 0.94137 & 391  & & 23 & 0.04960296200320839231 & 0.97809 &   \\[1ex]
9 & 0.07894768465249571895 & 0.94641 & 381  & & 24 & 0.04856391475086479136 & 0.97905 &   \\[1ex]
10 & 0.07506629658798923053 & 0.95084 & 371  & & 25 & 0.04758950454972984014 & 0.97994 &   \\[1ex]
11 & 0.07166627779626831642 & 0.95471 & 361  & & 26 & 0.04667334411591790177 & 0.98075 &   \\[1ex]
12 & 0.06866339082734667271 & 0.95810 & 351  & & 27 & 0.04580986165578441533 & 0.98150 &   \\[1ex]
13 & 0.06599134928348895227 & 0.96108 & 341  & & 28 & 0.04499417680286498705 & 0.98219 &   \\[1ex]
14 & 0.06359753741872661359 & 0.96373 &   & & 29 & 0.04422199837110163796 & 0.98284 &   \\[1ex]
15 & 0.06143976881171471065 & 0.96607 &   & & 30 & 0.04348953956755718180 & 0.98344 &   \\[1ex]
16 & 0.05948387383735141015 & 0.96817 &   & & 31 & 0.04279344727893972947 & 0.98399 &   \\[1ex]
\bottomrule
\end{tabular}
\end{table}

Second, at the other extreme where $\e\lesssim 1$, it becomes more convenient to
expand in powers of stretch factor, $\varepsilon=\sqrt{1-\e^2}$. 
The so-called \textit{asymptotic expansion}, valid for $\varepsilon\gtrsim 0$, is now,
\begin{multline}
A'=\pi\varepsilon:\qquad  \lambda'_0 = \sum_{\nu=-2}^\infty c_\nu \varepsilon^\nu = \frac{\pi^2}{4\varepsilon^2} +\frac{\pi}{2\varepsilon}+\frac{3}{4}
    + \left(\frac{11}{8\pi}+\frac{\pi}{12}\right)\, \varepsilon
\\[1ex]
  + \left(\frac{61}{16\pi^2}+\frac{1}{12}\right)\,\varepsilon^2
  + \left(\frac{1971}{128\pi^3}-\frac{9}{16\pi}+\frac{3\pi}{80}\right)\,\varepsilon^3
  + \left(\frac{20851}{256\pi^4}-\frac{271}{48\pi^2}+\frac{2}{45}\right)\,\varepsilon^4 
  \\[1ex]
  + \left(\frac{537219}{1024\pi^5}-\frac{11667}{256\pi^3}-\frac{7}{64\pi}+\frac{5\pi}{224}\right)\,\varepsilon^5
  + O(\varepsilon^6)\qquad
\label{eq:resultB}
\end{multline}
where, as anticipated, the first term corresponds to the infinite
strip eigenvalue per Eq.~(\ref{eq:limitstrip}).  The
coefficients~--~expressed as rational polynomials in powers of
$\pi$~--~are exact; but, for reference, rounded numerical values of the
first ten coefficients appear in Table~\ref{tab:resultBnumerical}.

This asymptotic expansion has a very different history. Troesch and
Troesch~\cite{tt1973,t1973} in 1973 appear to have started the
discussion\footnote{Although this ellipse problem has a very long
  history, they appear to be the first who promoted expanding the
  eigenvalue in powers of $\varepsilon$. In the second quoted article,
  they raise an intriguing Kac-like question: Can an elliptical drum
  have a spectrum similar to a stringed instrument?} and derived the
first three terms, up to zeroth order (i.e., 3/4) based on
\textit{expansions of the roots of the modified Mathieu functions for
  large eccentricity} (cf., \cite{tt1973} paragraph 6\footnote{In
  my notation, they had (with $m\!=\!k\!=\!0$),
  $(\lambda'_0)^{1/2}\varepsilon\!=\!\pi/2\!+\!\varepsilon/2\!+\!\varepsilon^2/(2\pi)\!+\!O(\varepsilon^3)$,
  and squaring, $\lambda'_0\varepsilon^2\!=\!\pi^2/4\!+\!\pi\varepsilon/2\!+\!3\varepsilon^2/4\!+\!O(\varepsilon^3)$.}).
Several decades later, in the mid-2000s, Borisov and
Freitas~\cite{bf2009} considered the ellipse as one example using a
boundary variation method. They published the first four terms, i.e.,
the first line of Eq.~(\ref{eq:resultB}), up to first order in
$\varepsilon$.

My second contribution to this problem is the next four terms of Eq.~(\ref{eq:resultB}),
i.e., the last two lines. 

\begin{table}[!htb]
\centering
\caption{\small Numerical values of the leading ten coefficients of the asymptotic
  series for the eigenvalue within the ellipse per
  Eq.~(\ref{eq:resultB}), 
  and rounded to twenty decimal
  places. Also indicated is $D_\nu$, the number of digits in agreement
  between the numerical coefficient (via these fits of data)
  and the respective coefficient displayed in Eq.~(\ref{eq:resultB}).  
  This truncated list of coefficients
  is based on polynomial interpolations using thirty-six 200-digit eigenvalues with eccentricity $\e=0.999800$
  to $0.999995$, each one incorporating the (known) closed-form lower-order coefficients.}
\label{tab:resultBnumerical}
\begin{tabular}{rlcccrlcc}
\toprule
$\nu$ &\qquad$c_\nu$ & $c_\nu/c_{\nu-1}$ & $D_\nu$ &\qquad \qquad & $\nu$ & \qquad$c_\nu$ & $c_\nu/c_{\nu-1}$\\[1ex]
\midrule
-2 & 2.46740110027233965471 & & 48 & & 3 & 0.43538365077995525294 & 0.92710 & 40 \\[1ex]
-1 & 1.57079632679489661923 & 0.63662 & 47 & & 4 & 0.30855816280914840552 & 0.70870 & 38 \\[1ex]
0 & 0.75000000000000000000 & 0.47746 & 45 & & 5 & 0.27983128169766678772 & 0.90690 & 36 \\[1ex]
1 & 0.69947548130186160990 & 0.93263 & 43 & & 6 & 0.19027912693622176700 & 0.67998 &  \\[1ex]
2 & 0.46962034596974608696 & 0.67139 & 41 & & 7 & 0.15981739762228202463 & 0.83991 &  \\[1ex]
\bottomrule
\end{tabular}
\end{table}

Of note is that as one moves down these series, the terms appear to
become more complicated, but intriguing patterns do
emerge~--~suggesting future work. For example, with the Maclaurin
series, the ratio of successive coefficients (see
Table~\ref{tab:resultAnumerical}) appears to approach unity (as
$\nu\!\to\!\infty$). This fact suggests that better a representation
is possible, quite analogous to the way that
$1\!+\!x\!+\!x^2\!+\!x^3\!+\!\cdots$ is better represented by
$1/(1\!-\!x)$. Despite that, the two canonical representations chosen
for this report are intended to help compare with prior results.

\section*{Technique}

Boundary variation techniques are mathematically elegant, but as
workers have discovered and acknowledged, getting simple results is
quite challenging. The other method~--~using Mathieu functions~--~is
also quite elegant, but it too leads to interesting computational
challenges.

The method I use is quite different from those other methods. I view
it more as a brute force method, however, it is guided and motivated
by those prior results.

My starting point is a high-precision computation of the fundamental
eigenvalue for judiciously chosen values of eccentricity.  For the
ellipse eigenvalues, I use the same method~\cite{j2017} as popularized
by Fox, Henrici, and Moler~\cite{fhm1967} fifty years earlier. The key
to make it work well is to use multi-precision software and pay
attention to the spacing of boundary points.\footnote{\raggedright
  Incidentally, the \textit{bounding method} of~\cite{j2017} does not
  work since the approximate eigenvalues do not alternate as more
  boundary matching points are added; however, the method is
  sufficiently convergent to yield multi-hundred-digit results in a
  relatively short time (e.g., 10,000 100-digit eigenvalues in a
  couple of days on a laptop for equally-spaced eccentricity from zero
  to 0.9999.). Since the ellipse boundary meets the axes in the first
  quadrant at right-angles, this provides a good clue as to why that
  bounding method does work for shapes where that doesn't happen.}

Next, sets of computed eigenvalues (``data points'') are fit to an
appropriate model equation and high-precision coefficients are
computed. By adding terms (and data points), the precision of each
coefficient can be improved, and thus estimated.\footnote{Knowing that
  a set of eigenvalues is good to 500 digits, effective eigenvalue
  bounds may be created by adding and subtracting a relative
  difference of $10^{-500}$ from those eigenvalues. Coefficients are
  then computed separately for the sets of upper and lower bounds,
  from which coefficient precision may be estimated.} It is quite
typical to work with coefficients that appear to be precise to more
than a few dozen digits, and often to hundreds of digits.

When the leading coefficient is of sufficient precision, it is fed
into an LLL integer relation algorithm with an ansatz guided by prior
results.  If the precision of the coefficient is too low, the LLL
output is far from unique. As the coefficient precision increases, a
viable candidate for the integer relation clearly emerges, and with
increasing confidence.  After numerical evidence supports a result for
a coefficient, it is incorporated into the dependent variable of the
model equation, and the process repeated in search of the next
coefficient.

A very similar method works quite well for the $1/S$ expansion of the
fundamental Dirichlet Laplacian eigenvalue of the $S$-sided regular
polygon~\cite{j2017a}.

To get the coefficients in the Maclaurin series,
Eq.~(\ref{eq:resultA}) with Table~\ref{tab:resultA}, thirty 500-digit
eigenvalues, equally spaced from $\e=0.000001$ to~$0.000030$, were
sufficient. Each eigenvalue required about five hours on my commodity
hardware (i7, 6-core, desktop) using free software (pari/gp). A
thirty-term interpolating polynomial in even powers of $\e$, with
those thirty eigenvalues, was then determined. Numerical values of the 
first thirty non-trivial coefficients (using sixty eigenvalues) of the
Maclaurin series are listed in Table~\ref{tab:resultAnumerical}.

To get the coefficients in the asymptotic expansion,
Eq~(\ref{eq:resultB}), some twenty 200-digit eigenvalues in the
interval from $\e=0.999800$ to $0.999995$ were sufficient. That entire eigenvalue
computation took a few days with the same setup.  Like the Maclaurin
series, an interpolating polynomial, this time in powers of
$\varepsilon$, was used.

With a sufficiently precise (numerical) coefficient in hand, the technique to
determine an expression is given by example, here for the third-order
term in Eq.~(\ref{eq:resultB}), i.e., the coefficient $c_3$
multiplying $\varepsilon^3$. 
Guided by lower order expressions, an ansatz might look like
\begin{equation}
  a_1\,c_3 + \frac{a_2}{\pi^5} + \frac{a_3}{\pi^3} + \frac{a_4}{\pi} + a_5\,\pi + a_6\,\pi^3
  = 0
  \end{equation}
where the six ``small'' integers $a_i$ are to be sought using LLL. In this
particular example, if the LLL output has $a_2=a_6=0$, that adds more
confidence. Fitting the numerical eigenvalue data to the model equation indicates that 
\begin{equation}
  c_3 = 0.\underline{435383650779955252940603845025457624}39389188\cdots
\end{equation}
where the underlined leading 36 digits appear to be
correct (based only on the fit to the data). With that, the LLL
routine (unambiguously) yields,
\begin{equation}
  a_1=640\qquad a_2=0 \qquad a_3=-9855\qquad a_4=360\qquad a_5=-24\qquad a_6=0
\end{equation}
from which the coefficient is constructed. Incidentally, this constructed coefficient reproduces the
numerical coefficient,
\begin{equation}
 \left(\frac{1971}{128\pi^3}-\frac{9}{16\pi}+\frac{3\pi}{80}\right) =
  0.\underline{435383650779955252940603845025457624}43837546\cdots
\end{equation}
matching the above $c_3$ to 36 digits. This then becomes the convincing evidence (valid to 36 digits)
that the LLL gave the correct result.

Two effects provide a practical limit to this method. Given a set of
computed eigenvalues, as one progresses along a series: (1)~numerical
precision of the coefficients decreases and (2)~the number of terms
(rational polynomial in powers of $\rho$ or $\pi$) needed to represent
a coefficient increases. At some point, these two effects conspire in
such a way that LLL routine simply fails to offer an unambiguous
solution. All of the results in this report were pushed to the limit
and only unambiguous, unique solutions are presented. The simple way to
extend results even further is with a more extensive eigenvalue computation.

\section*{Conclusion}

The method outlined in this report appears to be quite fruitful in
confirming and extending power series expansions of the eigenvalues of
the Laplacian for the ellipse.  The main results of this report are
summarized in Eq.~(\ref{eq:resultA}) with Table~\ref{tab:resultA}, and
Eq.~(\ref{eq:resultB}).  This method nicely complements other methods,
such as boundary variation methods, providing independent confirmation
and hints on what terms look like on down those series.

The ingredients necessary for the method to work are (1)~very high
precision eigenvalue computations, (2)~an analytical model to give
hints at what a series might look like (in terms of $\pi$ or
$\besselroot$, for example), (3)~an integer relation algorithm, and
(4)~a little patience and luck.

Of course, having these results raises many other questions,
especially as simple patterns are exposed. Questions might include:
What do higher eigenvalues look like? What are the convergence
properties? How do these results relate to Mathieu functions? Are
there better representations~--~such as rational polynomials or
continued fractions? What are the recurrence relations for the
coefficients?  Can these results help guide boundary variation
techniques that solve the same problem? Indeed, there seem
to be more new questions raised than answers provided.

\section*{Appendix}
Refer to Fig.~\ref{fig:ellipse} for artwork. To fix the notation,
\begin{equation}
  \left(\frac{x}{a}\right)^2+\left(\frac{y}{b}\right)^2=1\qquad
  \varepsilon=\frac{a}{b}=\sqrt{1-\e^2}\qquad
  \e=\sqrt{1-\varepsilon^2}\qquad 
  A=\pi ab \qquad
  \label{eq:ellipse}
\end{equation}
where $a$ and $b$ are, respectively, the semi-major and semi-minor
axes, $\varepsilon$ the stretch factor, $\e$ the eccentricity, and $A$
the area.  Without loss of generality, it is tacitly assumed that
$a\ge b$ so that as eccentricity varies from zero to unity, the stretch
factor varies from unity to zero.

Key to sorting out the expressions in this report~--~and comparing to
other results~--~is the invariant product (eigenvalue $\times$ area)
for a given eccentricity. I shall use two popular conventions, 
distinguished using a prime, and here showing which quantities depend on eccentricity,
\begin{subequations}
  \begin{align}
\text{Constant area:} \quad & A=\pi \qquad \lambda_0(\e) \qquad
a(\e)\cdot b(\e)=1\\ \text{Constant semi-major axis:}\quad & A'(\e)=\pi\varepsilon(\e)
\qquad \lambda'_0(\e)=\lambda_0(\e)/\varepsilon(\e)\qquad a'=1 \qquad
b'(\e)=\varepsilon(\e)
\end{align}
\end{subequations}
where, in \textit{both} cases, the unit-radius circle is the shape when
$\e=0$, and ``constant'' means independent of $\e$.
The invariant product $A\lambda_0=A'\lambda'_0$ means
$\lambda'_0\varepsilon=\lambda_0$.
The numerical values of both $\lambda_0$ and $\lambda'_0$ diverge as
$\e\to 1^{-}$, but the product $\lambda'_0\varepsilon^2=\lambda_0\varepsilon= (\pi/2)^2$ in
that limit, see Eq.~(\ref{eq:limitstrip}).  

To ameliorate the inevitable confusion, the ellipse area
shall always be clearly specified for each eigenvalue expression or
set of data.

\bibliographystyle{plain}
\bibliography{references}
\end{document}